\documentclass[12pt]{amsart}
\pdfoutput=1
\usepackage{graphicx}
\usepackage{amssymb}
\usepackage{amsmath}
\usepackage{amsthm}
\usepackage{amscd}
\usepackage{epsfig}
\usepackage{color}

\allowdisplaybreaks[4]
\newtheorem{theorem}{Theorem}[section]

\newtheorem{lemma}[theorem]{Lemma}

\newtheorem{proposition}[theorem]{Proposition}

\theoremstyle{definition}
\newtheorem{definition}[theorem]{Definition}
\numberwithin{equation}{section}

\newcommand{\N}{\mathbb N}

\begin{document}

\title{Packing entropy for fixed-point free flows}

\author [Ruiming Liang and Haoyi Lei ]{Ruiming Liang and Haoyi Lei }

\address{Department of Mathematics, Nanjing University,
Nanjing, Jiangsu, 210093, P.R. China} \email{181180069@smail.nju.edu.cn}

\address{Department of Mathematics, Nanjing University,
Nanjing, Jiangsu, 210093, P.R. China} \email{181870091@smail.nju.edu.cn}

\keywords {measure-theoretic entropy, packing entropy, variational principle, fixed-point free flow, reparametrization}

\begin{abstract}
Let $(X,\phi)$ be a compact flow without fixed points. We define the packing topological entropy $h_{\mathrm{top}}^P(\phi,K)$ on subsets of $X$ through considering all the possible reparametrizations of time. For fixed-point free flows, we prove the following result: for any non-empty compact subset $K$ of $X$, $$h_{\mathrm{top}}^P(\phi,K)=\sup\{\overline{h}_{\mu}(\phi):\mu(K)=1,\mu\text{ is a Borel probability measure on} X\},$$ where $\overline{h}_{\mu}(\phi)$ denotes the upper local entropy for a Borel probability measure $\mu$ on $X$.
\end{abstract}

\maketitle


\section{Introduction}

Entropy plays a crucial role in the study of the complexity of dynamical systems, which characterizes the growth of geometric or measure-theoretic quantities under long time dynamical iterations. Among various definitions of entropies, one of frequently studied definitions was introduce by Bowen \cite{B3} in 1973 in a way resembling the Hausdorff dimension for discrete dynamical systems, which is now called Bowen topological entropy. Bowen topological entropy can be defined on any subset of the state space while the classic topological entropy defined by Adler, Konheim and McAndrew can not, although they coincide on the whole state space. Similarly, by considering the concept of packing dimension, Feng and Huang introduced the concept of packing topological entropy over subsets for discrete systems \cite{FH}. In particular, it also coincides with Adler-Konheim-McAndrew's topological entropy on the whole space. This two definitions are of great importance in topological dynamics and dimension theory. Inspired by classical results in geometric-measure theory, Feng and Huang proved the following non-classic variational principles for Bowen entropy and packing entropy:

Let $(X,T)$ be a topological dynamical system where $X$ is a compact metric space and $T: X\rightarrow X$ is a continuous map. For any non-empty subset $K$ of $X$, let $h_{\mathrm{top}}^B(T,K)$ and $h_{\mathrm{top}}^P(T,K)$ denote the Bowen entropy and packing entropy on $K$, respectively.
When $K$ is a compact subset of $X$,
\begin{align}
 \label{eq-1} h_{\mathrm{top}}^B(T,K)=\sup\{\underline{h}_{\mu}(T):\mu(K)=1,\mu\in M(X)\};\\
 \label{eq-2} h_{\mathrm{top}}^P(T,K)=\sup\{\overline{h}_{\mu}(T):\mu(K)=1,\mu\in M(X)\}.
\end{align}
Here $M(X)$ denotes the set of all Borel probability measures on $X$ and $\underline{h}_{\mu}(T)$ and $\overline{h}_{\mu}(T)$ are the lower and upper local entropies for $\mu$ (see \cite{FH} for more details).

Feng and Huang's variational principles have now been extended to many situations such as for more general dynamical systems 
under countable amenable group actions \cite{ZC,DZZ} and for nonautonomous dynamical systems \cite{XZ}. A basic but nature question is how the above result works for the case of continuous dynamical systems (flows).
Let $(X,d)$ be a compact metric space with metric $d$. Recall that a pair $(X,\phi)$ is called a flow, if $\phi:X\times\mathbb{R}\to X$ is a continuous map satisfying $\phi_t\circ\phi_s=\phi_{s+t}$ for all $s,t\in\mathbb{R}$ and $\phi_{t}(\cdot)=\phi(\cdot,t)$ is a homeomorphism on $X$.
A Borel probability measure $\mu$ on $X$ is called {\it $\phi-$invariant} if for any Borel set $B$, it holds $\mu(\phi_{t}(B))=\mu(B)$ for all $t\in\mathbb{R}$. It is called {\it ergodic} if any $\phi-$invariant Borel set has measure 0 or 1. For convention, we denote all $\phi-$invariant Borel probability measures and all ergodic $\phi-$invariant Borel probability measures on $X$ by $M(X,\phi)$ and $E(X,\phi)$ respectively.

A nature definition of entropy for flows is to study the time-1 map. This deduces the entropy theory for flows to discrete case. But for flows there are more depth considerations--finding quantitative invariants for flows under orbit equivalents. This leads to a series of studies of flows through considering all the possible reparametrizations of flows (\cite{T1,T2,S,SV}). Recently Dou etc. in \cite{DFQ} introduced Bowen topological entropy for flows via the reparametrization balls and proved the variational principle \eqref{eq-1} for flows without fixed-points. For further studies on Bowen topological entropy for fixed-point free flows, see \cite{JC, WC}. In the present paper, we will study packing topological entropy for fixed-point free flows. Especially we will aim on proving the variational principle \eqref{eq-2} for flows without fixed-points.

The paper is organized as follows. In section 2, we introduce packing topological entropy for flows via reparametrization balls and some basic properties are also listed there. In section 3, we give the exact statement of our variational principle and some preparatory lemmas for the proof. Finally, in section 4, we give detailed proof of the variational principle of packing topological entropy for flows without fixed-points. For the proof of the lower bound, we need carefully treatments on the reparametrization balls. For the proof of the upper bound, we will employ Feng and Huang's method in \cite{FH}, the idea of which origins from Joyce and Preiss's work on packing measures \cite{JP}. This idea from geometric measure theory is surprisingly suitable for many cases in dynamical systems.

\section{Packing topological entropy via reparametrization balls}

To introduce the concept of packing topological entropy for subsets of a flow, we first need to introduce the definition of a reparametrization ball, which is an analogy of the traditional Bowen ball.

\begin{definition}

For a closed interval $I$ which contains the origin, a continuous map $\alpha: I\rightarrow \mathbb{R}$ is called a {\it reparametrization} if it is a homeomorphism onto its image and $\alpha(0)=0$. The set of all such reparametrizations on $I$ is denoted by $Rep(I)$. For a flow $\phi$ on $X$, $x\in X$, $t\in \mathbb{R}^+$ and $\varepsilon>0$, we set
\begin{align*}
B(x,t,\varepsilon,\phi)=\{y\in X: \text{ there exists } \alpha \in Rep[0,t] \text{ such that } \\
                    d(\phi_{\alpha(s)}{x},\phi_{s}y)<\varepsilon, \text{ for all }0\le s\le t\},
\end{align*}

and
\begin{align*}
\overline{B}(x,t,\varepsilon,\phi)=\{y\in X: \text{ there exists } \alpha \in Rep[0,t] \text{ such that } \\
                    d(\phi_{\alpha(s)}{x},\phi_{s}y)\leq\varepsilon, \text{ for all }0\le s\le t\},
\end{align*}

Usually, we call $B(x,t,\varepsilon,\phi)$ a {\it $(t,\varepsilon,\phi)-$ball} or a {\it reparametrization ball} in $X$. Clearly, all the reparametrization balls are open sets.

\end{definition}

Now, we can give the definition of packing topological entropy defined through reparametrization balls.

\begin{definition}

Let $(X,\phi)$ be a flow and $Z$ a subset of $X$. For $s\geq 0$, $N\in \mathbb{N}$, and $\varepsilon > 0$, define
\begin{displaymath}
P_{N,\varepsilon}^{s}(\phi,Z)=\sup \sum_{i}\exp(-st_{i}),
\end{displaymath}
where the supremum is taken over all finite or countable families of disjoint closed reparametrization balls $\{\overline{B}(x_{i},t_{i},\varepsilon,\phi)\}$ such that $x_{i}\in X$, $t_{i}\ge N$.

The quantity $P_{N,\varepsilon}^{s}$ dose not increase as $N$ increases, hence the following limit exists:
$$P_{\varepsilon}^{s}(\phi,Z)=\lim_{N\rightarrow \infty}P_{N,\varepsilon}^{s}(\phi,Z).$$

Define $\mathcal{P}_{\varepsilon}^s(Z)=\inf\{\sum\limits_{i=1}^\infty P_{\varepsilon}^{s}(\phi,Z_i): Z\subset\bigcup\limits_{i=1}^\infty Z_i\}$.

By definition and simple calculation, we know $$h_{\mathrm{top}}^P(\phi,Z,\varepsilon)=\inf\{s:\mathcal{P}_{\varepsilon}^s(Z)=0\}=\sup\{s:\mathcal{P}_{\varepsilon}^s(Z)=+\infty\}$$ is well-defined.

Note that $h_{\mathrm{top}}^P(\phi,Z,\varepsilon)$ does not decrease as $\varepsilon$ decreases. Define the {\it packing topological entropy} for subset $Z$ of $X$ to be $$h_{\mathrm{top}}^P(\phi,Z)=\lim\limits_{\varepsilon\to0}h_{\mathrm{top}}^P(\phi,Z,\varepsilon).$$
\end{definition}

For packing topological entropy we have the following properties.

\begin{proposition}
Let $Z,Z'$ and $Z_1,Z_2,\ldots$ be subsets of $X$.
\begin{enumerate}
  \item If $Z\subset Z'$, then $h_{\mathrm{top}}^P(\phi,Z)\leq h_{\mathrm{top}}^P(\phi,Z').$
  \item If $Z\subset\bigcup\limits_{i=1}^\infty Z_i$,  then for $s\geq 0$ and $\varepsilon>0$, we have $\mathcal{P}_\varepsilon^s(Z)\leq\sum\limits_{i=1}^\infty\mathcal{P}_\varepsilon^s(Z_i)$ and
      $h_{\mathrm{top}}^P(\phi,Z)\leq\sup\limits_{i\geq1}h_{\mathrm{top}}^P(\phi,Z_i).$
\end{enumerate}
\end{proposition}
The proof is simple from the definition and we omit it.

Since the Hausdorff dimension and the packing dimension are dual concepts in fractal geometry, we will compare packing topological entropy with Bowen topological entropy for a flow. In the following we recall the definition of Bowen topological entropy for a flow (see \cite{DFQ}).

\begin{definition}
Let $(X,\phi)$ be a flow and $Z$ a subset of $X$. For $s\ge 0$, $N\in \mathbb{N}$ and $\varepsilon > 0$, define
\begin{displaymath}
\mathcal{M}_{N,\varepsilon}^{s}(\phi,Z)=\inf \sum_{i}\exp(-st_{i}),
\end{displaymath}
where the infimum is taken over all finite or countable families of reparametrization balls $\{B(x_{i},t_{i},\varepsilon,\phi)\}$ such that $x_{i}\in X$,  $t_{i}\ge N$ and
$\bigcup B(x_{i},t_{i},\varepsilon,\phi)\supset Z$.

The quantity $\mathcal{M}_{N,\varepsilon}^{s}$ dose not decrease as $N$ increases and $\varepsilon$ decreases, hence the following limits exist:
$$\mathcal{M}_{\varepsilon}^{s}(\phi,Z)=\lim_{N\rightarrow \infty}\mathcal{M}_{N,\varepsilon}^{s}(\phi,Z),
\mathcal{M}^{s}(\phi,Z)=\lim_{\varepsilon\rightarrow 0}\mathcal{M}_{\varepsilon}^{s}(\phi,Z).$$

The {\it Bowen topological entropy} $h_{top}^{B}(\phi,Z)$ is defined as a critical value of the parameter $s$, where $\mathcal{M}^{s}(\phi,Z)$ jumps from $\infty$ to $0$, i.e.
\begin{align*}
h_{top}^{B}(\phi,Z)&=\inf \{s:\mathcal{M}^{s}(\phi,Z)=0\}\\
           &=\sup \{s:\mathcal{M}^{s}(\phi,Z)=\infty\}.
\end{align*}
\end{definition}

We now list some properties of Bowen topological entropy (\cite[Proposition 2]{DFQ}).
\begin{proposition}
Let $Z,Z'$ and $Z_1,Z_2,\ldots$ be subsets of $X$.
\begin{enumerate}
  \item If $Z\subset Z'$, then $h_{\mathrm{top}}^B(\phi,Z)\leq h_{\mathrm{top}}^B\phi,Z').$
  \item If $Z\subset\bigcup\limits_{i=1}^\infty Z_i$,  then for $s\geq 0$ and $\varepsilon>0$, we have $\mathcal{M}_\varepsilon^s(Z)\leq\sum\limits_{i=1}^\infty\mathcal{M}_\varepsilon^s(Z_i)$ and
      $h_{\mathrm{top}}^B(\phi,Z)\leq\sup\limits_{i\geq1}h_{\mathrm{top}}^B(\phi,Z_i).$
\end{enumerate}
\end{proposition}

Moreover, we have
\begin{proposition}
  For any $Z\subset X, h_{\mathrm{top}}^B(\phi,Z)\leq h_{\mathrm{top}}^P(\phi,Z)$.
\end{proposition}
\begin{proof}
  If $h_{\mathrm{top}}^B(\phi,Z)=0$ then there is nothing to prove. Assume $h_{\mathrm{top}}^B(\phi,Z)>0$ and $0<s<h_{\mathrm{top}}^B(\phi,Z)$. For any $n\in\N$ and $\varepsilon>0$,
  let $\{\overline{B}(x_{i},n,\varepsilon,\phi)\}_{i=1}^R$ be a disjoint family with $x_i\in Z$ such that the cardinality $R=R_n(Z,\varepsilon)$
  is maximal. Then for any $\delta>0$, $\bigcup_{i=1}^R\overline{B}(x_{i},n,2\varepsilon+\delta,\phi)\supset Z$. Hence
  $$\mathcal{M}_{n,2\varepsilon+\delta}^{s}(\phi,Z)\le R\exp(-sn)\le P_{n,\varepsilon}^s(\phi, Z).$$
  Letting $n\rightarrow \infty$, we have $\mathcal{M}_{2\varepsilon+\delta}^{s}(\phi,Z)\le P_{\varepsilon}^s(\phi, Z)$.
  Thus for any $\cup_{i=1}^{\infty}Z_i\supset Z$,
  $$\mathcal{M}_{2\varepsilon+\delta}^{s}(\phi,Z)\le \sum_{i=1}^{\infty}\mathcal{M}_{2\varepsilon+\delta}^{s}(\phi,Z_i)\le\sum_{i=1}^{\infty} \mathcal{P}_{\varepsilon}^s(\phi, Z_i),$$
  which implies that $\mathcal{M}_{2\varepsilon+\delta}^{s}(\phi,Z)\le \mathcal{P}_{\varepsilon}^s(\phi, Z)$.
  Note that $0<s<h_{\mathrm{top}}^B(\phi,Z)$. Hence $\mathcal{M}^s(\phi, Z)=\infty$ and $\mathcal{M}_{2\varepsilon+\delta}^{s}(\phi,Z)>1$ when
  $\varepsilon$ and $\delta$ are sufficiently small. Hence $\mathcal{P}_{\varepsilon}^s(\phi, Z)>1$ and $h_{\mathrm{top}}^P(\phi,Z,\varepsilon)\ge s$
  when $\varepsilon$ is small. Since $h_{\mathrm{top}}^P(\phi,Z,\varepsilon)$ does not decrease as $\varepsilon$ decreases,
  $h_{\mathrm{top}}^P(\phi,Z)=\lim\limits_{\varepsilon\to 0}h_{\mathrm{top}}^P(\phi,Z,\varepsilon)\ge s$. Therefor $h_{\mathrm{top}}^P(\phi,Z)\ge h_{\mathrm{top}}^B(\phi,Z)$.
\end{proof}

\if
Moreover, we have
\begin{proposition}
\begin{enumerate}
  \item For any $Z\subset X, h_{\mathrm{top}}^B(\phi,Z)\leq h_{\mathrm{top}}^P(\phi,Z)$.
  \item If $Z$ is $\phi-$invariant and compact, then $$h_{\mathrm{top}}^B(\phi,Z)=h_{\mathrm{top}}^P(\phi,Z).$$
  \item Let $h_{\mathrm{top}}(\phi,X)$ denote the traditional topological entropy defined by time-1 mapping. If $(X,\phi)$ is a compact flow without fixed point, then we have $$h_{\mathrm{top}}^B(\phi,X)=h_{\mathrm{top}}^P(\phi,X)=h_{\mathrm{top}}(\phi,X).$$
\end{enumerate}
\end{proposition}
\fi

\section{Statement of main theorem and preparatory lemmas}

The measure-theoretic local entropies are defined as follows.
\begin{definition}
Let $\mu\in \mathcal{M}(X)$. The {\it measure-theoretic lower and upper local entropies} of $\mu$ are defined respectively by
\begin{equation*}
\underline{h}_{\mu}(\phi)=\int\underline{h}_{\mu}(\phi,x)\,d\mu,\text{ and }
\overline{h}_{\mu}(\phi)=\int\overline{h}_{\mu}(\phi,x)\,d\mu
\end{equation*}
where
$$\underline{h}_{\mu}(\phi,x)=\lim_{\varepsilon\rightarrow0}\liminf_{t\rightarrow +\infty}-\frac{1}{t}\log\mu(B(x,t,\varepsilon,\phi))$$
and
$$\overline{h}_{\mu}(\phi,x)=\lim_{\varepsilon\rightarrow0}\limsup_{t\rightarrow +\infty}-\frac{1}{t}\log\mu(B(x,t,\varepsilon,\phi)).$$

\end{definition}

Now we state the main theorem.

\begin{theorem} \label{main}
Let $(X,\phi)$ be a compact metric flow without fixed points. If $K$ is a non-empty compact subset of $X$, then
$$
h_{top}^{P}(\phi,K)=\sup\{\overline{h}_{\mu}(\phi):\mu\in M(X),\mu(K)=1\}.
$$
\end{theorem}

We suggest here that there are some further results related to Theorem \ref{main} for flows without fixed points. Due to Feng and Huang \cite{FH}, the compact subsets $K$'s can be improved to any analytic subset of $X$ when the classic topological entropy of the flow is finite. (We note that here for the flows without fixed points, the classic topological entropy coincides with the topological entropy defined through reparametrization balls \cite{T1,T2}.)

In this section, we first will give some properties about reparametrization balls for flows without fixed points and then give a covering lemma. These lemmas are crutial in proving
Theorem \ref{main}.

\begin{lemma}[Lemma 1.2 of \cite{T1}] \label{lemma1}
Let $(X,\phi)$ be a compact metric flow without fixed points. For any $\eta>0$, there exists $\theta>0$ such that for any
$x,y\in X $, any interval $I$ containing the origin, and any reparametrization $\alpha\in Rep(I)$, if $d(\phi_{\alpha(s)}(x),\phi_{s}(y))<\theta$ for all $s\in I$, then it holds that
\begin{equation*}
|\alpha(s)-s|<\begin{cases} \eta|s|, &\text{ if }|s|>1, \\ \eta, &\text{ if } |s|\le 1.
\end{cases}
\end{equation*}
\end{lemma}

\if
Let $(X,\phi)$ be a compact metric flow without fixed points.
For any $\varepsilon>0$, there exists $\delta>0$ depending only on $\varepsilon$, such that
\begin{equation}
B(x,t_1,\varepsilon,\phi)\subset B(x,t_2,2\varepsilon,\phi), \text{ for any }x\in X,
\end{equation}\label{eq7}
whenever $t_1,t_2>0$ and $|t_1-t_2|<\delta$.
\fi

The following $5r$-lemma for reparametrization balls is proved by Dou etc. \cite[Theorem 3.5]{DFQ}. This is a variation of the classic $5r$-covering lemma and it will
play a crucial role for the proof of Theorem \ref{main}.

\begin{theorem}[{\bf $5r$-lemma for reparametrization balls}]\label{5r}
Let $(X,\phi)$ be a compact metric flow without fixed points. For $0<\eta<1$, let $\theta>0$ be as in Lemma \ref{lemma1}.
Let $\mathcal{B}=\{B(x,t,\varepsilon,\phi)\}_{(x,t)\in \mathcal{I}}$ be a family of reparametrization balls in X with $0<\varepsilon<\frac{\theta}{2}$ and $t>\frac{1}{(1-\eta)^2}$.
Then there exists a finite or countable subfamily $\mathcal{B}'=\{B(x,t,\varepsilon,\phi)\}_{(x,t)\in \mathcal {I}'}$($\mathcal{I}'\subset \mathcal{I}$) of pairwise disjoint reparametrization balls in $\mathcal{B}$ such that
\begin{equation*}
\underset{B\in\mathcal{B}}{\bigcup}B\subseteq\underset{(x,t)\in\mathcal{I}'}{\bigcup}B(x,\hat{t},5\varepsilon,\phi)
\end{equation*}
where 
$\hat{t}=(1-\eta)^2t$.
\end{theorem}

At the end of this section, we give a frequently used lemma when proving Theorem \ref{main}. It is analogous to \cite[Lemma 4.1]{FH}.

\begin{lemma}\label{lemma6}
Let $Z\subset X$ and $s,\varepsilon>0$. Assume $P_{\varepsilon}^s(Z)=\infty$. Then for any given finite interval $(a,b)\subset\mathbb{R}$ with $a\geq0$ and any $N\in\mathbb{N}$, there exists a finite  disjoint collection $\{\overline{B}(x_{i},t_{i},\varepsilon,\phi)\}$ such that $x_i\in Z, t_i\geq N$ and $\sum\limits_i e^{-t_i s}\in(a,b).$
\end{lemma}
\begin{proof}
The proof is the same with that of \cite[Lemma 4.1]{FH} for discrete dynamical systems. For the completeness we give the proof.

Let $N_1>N$ be sufficiently large such that $e^{-N_1s}<b-a$. Since $P_{N,\varepsilon}^{s}(\phi,Z)$ does not increase as $N$ increases and $P_{\varepsilon}^{s}(\phi,Z)=\lim\limits_{N\rightarrow \infty}P_{N,\varepsilon}^{s}(\phi,Z)$, we have $P_{N_1,\varepsilon}^s(Z)\ge P_{\varepsilon}^s(Z)=\infty.$ Thus we can choose a finite disjoint collection $\{\overline{B}(x_{i},t_{i},\varepsilon,\phi)\}$ with $x_i\in Z, t_i\geq N_1$ such that $\sum\limits_i e^{-t_is}>b$. Since $e^{-t_is}<b-a$, we can discard elements in this collection one by one until we have  $\sum\limits_i e^{-t_i s}\in(a,b).$
\end{proof}

\section{Proof of Theorem \ref{main}}

With the preparation in Section 3, we can now prove Theorem \ref{main}. The following lemma gives the proof of the lower bound.

\begin{lemma}\label{lemma4-1}
Let $(X,\phi)$ be a compact metric flow without fixed points. We have $$h_{top}^{P}(\phi,Z)\geq \sup\{\overline{h}_{\mu}(\phi):\mu\in\mathcal{M}(X),\mu(Z)=1\}$$ for any Borel set $Z\subset X$.
\end{lemma}

\begin{proof}

Let $\mu\in M(X)$ with $\mu(Z)=1$ for some Borel set $Z\subset X$. We just need to prove $h_{top}^{P}(\phi,Z)\geq \overline{h}_{\mu}(\phi)$. Since the lemma automatically holds when $\overline{h}_{\mu}(\phi)=0$, we may assume $\overline{h}_{\mu}(\phi)>0$. Let $0<s<\overline{h}_{\mu}(\phi)$, we shall prove $h_{top}^{P}(\phi,Z)\geq s.$

Choose $0<\eta<1$ such that $\frac{1}{(1-\eta)^2}<2$ and let $\theta>0$ be given as in Lemma \ref{lemma1} and Lemma \ref{5r}.

Note that $\mu(B(x,t,\varepsilon,\phi))$ does not increase as $\varepsilon$ decreases. Then there exists $0<\varepsilon<\frac{\theta}{2}$, $\delta>0$ and a Borel set $A\subset Z$ with $\mu(A)>0$ such that $$\overline{h}_\mu(\phi,x,\varepsilon)>s+\delta,\;\forall x\in A,$$
where $\overline{h}_\mu(\phi,x,\varepsilon)=\varlimsup\limits_{t\to\infty}-\frac{1}{t}\log\mu({B}(x,t,\varepsilon,\phi))$.
To see this, we note that $\overline{h}_{\mu}(\phi)=\int\overline{h}_{\mu}(\phi,x)\,d\mu>s$. And hence there exists $\delta>0$ and a Borel set $A$ with $\mu(A)>0$ such that for any $x\in A$, it holds that $\overline{h}_\mu(\phi,x)>s+\delta$. Since $\overline{h}_\mu(\phi,x)=\lim\limits_{\varepsilon\to0}\overline{h}_\mu(\phi,x,\varepsilon)$, there must exist a sufficiently small $\varepsilon>0$ which meets our requirement.

Next we show that $\mathcal{P}_{\frac{\varepsilon}{6}}^s(Z)=\infty,$ which implies that $h_{\mathrm{top}}^P(\phi,Z)\geq h_{\mathrm{top}}^P(\phi,Z,\frac{\varepsilon}{6})\geq s.$ It suffices to prove that $P_{\frac{\varepsilon}{6}}^s(E)=\infty$ for any Borel set $E\subset A$ with $\mu(E)>0$. Fix such a set $E$. Define $$E_t=\{x\in E:\mu(B(x,t,\varepsilon,\phi))<e^{-t(s+\delta)}\},\; t\in\mathbb{R}_{+},$$ 
and
\begin{align}\label{eq-En}
  \widetilde{E}_n=\{x\in E:\mu(B(x,n,\varepsilon,\phi))<e^{(-n+1)(s+\delta)})\},\;n\in\mathbb{N}.
\end{align}

Note that $B(x,\lceil t \rceil,\varepsilon,\phi)\subset B(x,t,\varepsilon,\phi)$. Hence for $x\in E_t$,
 $$\mu(B(x,\lceil t \rceil,\varepsilon,\phi))\leq\mu(B(x,t,\varepsilon,\phi))< e^{-t(s+\delta)}\leq e^{-\lceil t \rceil(s+\delta)}\cdot e^{s+\delta}.$$
 This implies $x\in\widetilde{E}_n$. And thus $E_t\subset \widetilde{E}_{\lceil t \rceil}.$

By $E\subset A,$ we have $E=\bigcup\limits_{t\geq N}E_t\subset\bigcup\limits_{n\geq N}\widetilde{E}_n\subset E $ for each $N\in\mathbb{N}$. Hence $\mu(\bigcup\limits_{n\geq N}\widetilde{E}_n)\geq\mu(E)$. Then there exists $n \geq N$ such that 
\begin{align}\label{eq-E}
  \mu(\widetilde{E}_n)\geq\frac{1}{n(n+1)}\mu(E).
\end{align}

Now we consider the family $\{B(x,\frac{n}{(1-\eta)^2},\frac{\varepsilon}{5},\phi):x\in \widetilde{E}_n\}$. By Lemma \ref{5r}, the $5$r-lemma for reparametrization balls, there exists a finite or countable union of disjoint family $\{B(x,\frac{n}{(1-\eta)^2},\frac{\varepsilon}{5},\phi)\}_{x\in\Lambda}$(where $\Lambda\subset \widetilde{E}_n$ is a finite or countable set) such that
\begin{align}\label{eq-5r}
  \widetilde{E}_n\subset\bigcup_{x\in \widetilde{E}_n}B(x,\frac{n}{(1-\eta)^2},\frac{\varepsilon}{5},\phi)\subset\bigcup_{x\in\Lambda}B(x,n,\varepsilon,\phi).
\end{align}

Hence
\begin{align*}
  P_{N,\frac{\varepsilon}{6}}^s(E)&\ge P_{N,\frac{\varepsilon}{6}}^s(\widetilde{E}_n)\ge\sum\limits_{{x\in\Lambda}}e^{-ns}\\
  &\qquad\qquad(\text{since }\{\overline B(x,\frac{n}{(1-\eta)^2},\frac{\varepsilon}{6},\phi)\}_{x\in\Lambda} \text{ is a disjoint family})\\
  &=e^{n\delta}\cdot e^{-(s+\delta)}\cdot\sum_{{x\in\Lambda}}e^{(-n+1)(s+\delta)}\\
  &\ge e^{n\delta-s-\delta}\cdot\sum_{{x\in\Lambda}}\mu(B(x,n,\varepsilon,\phi))\text{ (by \eqref{eq-5r})}\\
  &\ge e^{n\delta-s-\delta}\mu(\widetilde{E}_{n}) \text{ (by \eqref{eq-En})}\\
  &\ge \frac{e^{n\delta-s-\delta}}{n(n+1)}\mu(E)\text{ (by \eqref{eq-E})}.
\end{align*}
Let $N\to\infty$, then we have $P_{\frac{\varepsilon}{6}}^s(Z)=\infty$. This finishes the proof of the lemma.
\end{proof}

In the following we deal with the proof of the upper bound, which is more complicated compared with the proof of the lower bound. The core of the proof is to construct a sequence of positive finite measures supported on finite sets satisfying suitable conditions on reparametrization balls. 

\begin{lemma}\label{lemma4-2}
Let $(X,\phi)$ be a compact metric flow without fixed points. For any compact $K\subset X$ with $h_{top}^{P}(\phi,K)>0$ and any $0<s<h_{top}^{P}(\phi,K),$ there exists $\mu\in M(K)$ such that $\overline{h}_\mu(\phi)\geq s.$
\end{lemma}
\begin{proof}
Let $\varepsilon>0$ be sufficiently small such that $0<s<h_{top}^{P}(\phi,K,\varepsilon)$ and let $s< t< h_{top}^{P}(\phi,K,\varepsilon)$. 
We will inductively construct the following: 
\begin{enumerate}
  \item A sequence of finite sets $\{K_i\}_{i=1}^\infty$ with $K_i\subset K$.
  \item A sequence of finite measures $\{\mu_i\}_{i=1}^\infty$ with each $\mu_i$ is supported on $K_i$.
  \item A sequence of positive real-valued functions $\{m_i:K_i\to\mathbb{R}_+\}_{i=1}^{+\infty}$.
 \item A sequence of positive numbers $\{\gamma_i\}_{i=1}^{+\infty}$.
\end{enumerate}

Step 1. Construct $K_1$, $\mu_1$, $m_1(\cdot)$ and $\gamma_1$.

Note that $\mathcal{P}_{\varepsilon}^t(K)=\infty$. Let $$H=\bigcup\{G\subset X: G \text{ is open and }\mathcal{P}_{\varepsilon}^t(K\cap G)=0\}.$$
$H$ is an open subset of $X$ and it can be covered by countably many open sets $G$'s with $\mathcal{P}_{\varepsilon}^t(K\cap G)=0$. 
Then by the subadditivity of $\mathcal{P}_{\varepsilon}^t(\cdot)$, we immediately get $\mathcal{P}_{\varepsilon}^t(K\cap H)=0$.

Let $K'=K\backslash H=K\cap(X\backslash H)$. For an open set $G\subset X$, if $\mathcal{P}_{\varepsilon}^t(K'\cap G)=0$, 
then $\mathcal{P}_{\varepsilon}^t(K\cap G)\le\mathcal{P}_{\varepsilon}^t(G\cap K')+\mathcal{P}_{\varepsilon}^t(K\cap H)=0$. 
This implies that $G\subset H$. Hence we have that for any open set $G$, either $K'\cap G=\emptyset$ or $\mathcal{P}_{\varepsilon}^t(K'\cap G)>0$.

Since $\mathcal{P}_{\varepsilon}^t(K)\leq\mathcal{P}_{\varepsilon}^t(K')+\mathcal{P}_{\varepsilon}^t(K\cap H)=\mathcal{P}_{\varepsilon}^t(K')$, 
we have $\mathcal{P}_{\varepsilon}^t(K')=\mathcal{P}_{\varepsilon}^t(K)=\infty$. Then it follows that $\mathcal{P}_{\varepsilon}^s(K')=\infty$ and thus  $P_{\varepsilon}^s(K')=\infty$. By Lemma \ref{lemma6}, we can find a finite  set $K_1\subset K'$ and a positive real-valued function $m_1(x)$ on $K_1$ such that the collection $\{\overline{B}(x,\varepsilon,m_1(x),\phi)\}_{x\in K_1}$ is disjoint and $$\sum\limits_{x\in K_1}e^{-m_1(x)s}\in(1,2).$$

Define $\mu_1=\sum\limits_{x\in K_1}e^{-m_1(x)s}\delta_x$, where $\delta_x$ denotes the Dirac measure at $x$. Take a small $\gamma_1>0$ such that for any function $z:K_1\to X$ with $d(x,z(x))<\gamma_1$ for each $x\in K_1$, we have for each $x\in K_1,$

\begin{align}
&\bigg(\overline{B}(z(x),\gamma_1)\cup\overline{B}(z(x),m_1(x),\varepsilon,\phi)\bigg)\bigcap\nonumber \\
&\qquad\qquad\qquad\bigg(\bigcup\limits_{y\in K_1\backslash\{x\}}\overline{B}(z(y),\gamma_1)\cup\overline{B}(z(y),m_1(x),\varepsilon,\phi)\bigg)=\emptyset.\nonumber
\end{align}
Here and afterwards, $\overline{B}(x,\varepsilon)$ denotes the closed ball $\{y\in X: d(x,y)\leq\varepsilon\}.$
 
Step 2. Construct $K_2$, $\mu_2$, $m_2(\cdot)$ and $\gamma_2$.

For each $x\in K_1$, we will construct as in Step 1, a finite set $$E_2(x)\subset K\cap B(x,\frac{\gamma_1}{4})$$ and a real-valued function $$m_2:E_2(x)\to[\max\{m_1(y):y\in K_1\},\infty)$$ such that the elements in $\{\overline{B}(y,m_2(y),\varepsilon,\phi)\}_{y\in E_2(x)}$ are disjoint, and
$$\mu_1(\{x\})<\sum\limits_{y\in E_2(x)}e^{-m_2(y)s}<(1+2^{-2})\mu_1(\{x\}).$$

Fix $x\in K_1$ and set $F=K\cap B(x,\frac{\gamma_1}{4})$. Let 
$$H_x=\bigcup\{G\subset X: G \text{ is open and }\mathcal{P}_{\varepsilon}^t(F\cap G)=0\}.$$ 
Let $F'=F\backslash  H_x$. Then as in Step 1, we can show that $\mathcal{P}_{\varepsilon}^t(F')=\mathcal{P}_{\varepsilon}^t(F)>0$. Furthermore, $\mathcal{P}_{\varepsilon}^t(F'\cap G)>0$ for any open set $G$ with $G\cap F'\neq \emptyset$. Since $s<t$, $\mathcal{P}_{\varepsilon}^s(F')=\infty$. 
Hence by Lemma \ref{lemma6}, we can find a finite set $E_2(x)\subset F'$ and a map $m_2:E_2(x)\to[\max\{m_1(y):y\in K_1\},\infty)$ such that 
\begin{enumerate}
  \item the elements in $\{\overline{B}(y,m_2(y),\varepsilon,\phi)\}_{y\in E_2(x)}$ are disjoint;
  \item $\mu_1(\{x\})<\sum\limits_{y\in E_2(x)}e^{-m_2(y)s}<(1+2^{-2})\mu_1(\{x\}).$
\end{enumerate}

Since the family $\{\overline{B}(x,\gamma_1)\}_{x\in K_1}$ is disjoint, $E_2(x)\cap E_2(x')=\emptyset$ for different $x,x'\in K_1$.
Define
$$K_2=\bigcup\limits_{x\in K_1}E_2(x)\text{ and }\mu_2=\sum\limits_{y\in K_2}e^{-m_2(y)s}\delta_y.$$
The elements in $\{\overline{B}(y,m_2(y),\varepsilon,\phi)\}_{y\in K_2}$ are pairwise disjoint because the elements in every $\{\overline{B}(y,m_2(y),\varepsilon,\phi)\}_{y\in E_2(x)}$ are disjoint and elements from two different $E_2(x)$'s are also disjoint. By similar argument as in Step 1, we can take $0<\gamma_2<\frac{\gamma_1}{4}$ small enough such that for any function $z:K_2\to X$ with $d(x,z(x))< \gamma_2$ for each $x\in K_2,$ we have
\begin{align}
&\bigg(\overline{B}(z(x),\gamma_2)\cup\overline{B}(z(x),m_2(x),\varepsilon,\phi)\bigg)\bigcap\nonumber\\
&\qquad\qquad\qquad\bigg(\bigcup\limits_{y\in K_2\backslash\{x\}}\overline{B}(z(y),\gamma_2)\cup\overline{B}(z(y),m_2(y),\varepsilon,\phi)\bigg)=\emptyset\nonumber
\end{align}
for each $x\in K_2$.

Step 3. Assume that $K_i,\mu_i,m_i(\cdot)$ and $\gamma_i$ have been constructed for $i=1,2,...,p$. In particular, assume that for any function $z:K_p\to X$ with $d(x,z(x))<\gamma_p$ for each $x\in K_p$, we have
\begin{align}
&\bigg(\overline{B}(z(x),\gamma_p)\cup\overline{B}(z(x),m_p(x),\varepsilon,\phi)\bigg)\bigcap\nonumber\\
&\qquad\qquad\qquad\bigg(\bigcup\limits_{y\in K_p\backslash\{x\}}\overline{B}(z(y),\gamma_p)\cup\overline{B}(z(y),m_p(y),\varepsilon,\phi)\bigg)=\emptyset\nonumber
\end{align}
for each $x\in K_p.$ We will construct $K_i,\mu_i,m_i(\cdot)$ and $\gamma_i$ for $i=p+1$ in a way similar to Step 2.

Note that the elements in $\{\overline{B}(x,\gamma_p)\}_{x\in K_p}$ are pairwise disjoint. For each $x\in K_p$, since $\mathcal{P}_{\varepsilon}^t(K\cap B(x,\gamma_p))>0,$ we can construct as in Step 2, a finite set 
$$E_{p+1}\subset K\cap B(x,\frac{\gamma_p}{4})$$ and a positive real-valued function $$m_{p+1}:E_{p+1}(x)\to[\max\{m_p(y):y\in K_p\},\infty)$$ such that the elements in $\{\overline{B}(y,m_{p+1}(y),\varepsilon,\phi)\}_{y\in E_{p+1}(x)}$ are disjoint, and $$\mu_p(\{x\})<\sum\limits_{y\in E_{p+1}(x)}e^{-m_{p+1}(y)s}<(1+2^{-p-1})\mu_p(\{x\}).$$

Clearly $E_{p+1}(x)\cap E_{p+1}(x')=\emptyset$ for different $x,x'\in K_p.$ 
Define $$K_{p+1}=\bigcup\limits_{x\in K_p}E_{p+1}(x)\text{ and } \mu_{p+1}=\sum\limits_{y\in K_{p+1}}e^{-m_{p+1}(y)s}\delta_y.$$

Since the elements in $\{\overline{B}(y,m_p(y),\varepsilon,\phi)\}_{y\in K_p}$ are pairwise disjoint,  
we can take $0<\gamma_{p+1}<\frac{\gamma_p}{4}$ small enough such that for any function $z:K_{p+1}\to X$ with $d(x,z(x))<\gamma_{p+1}$ for $x\in K_{p+1},$ we have for each $x\in K_{p+1}$,

\begin{align}
&\bigg(\overline{B}(z(x),\gamma_{p+1})\cup\overline{B}(z(x),m_{p+1}(x),\varepsilon,\phi)\bigg)\bigcap\nonumber\\
&\qquad\qquad\qquad\bigg(\bigcup\limits_{y\in K_{p+1}\backslash\{x\}}\overline{B}(z(y),\gamma_{p+1})\cup\overline{B}(z(y),m_{p+1}(y),\varepsilon,\phi)\bigg)=\emptyset.\nonumber
\end{align}

As in the above steps, we can construct inductively the sequences $\{K_i\},\{\mu_i\},\{m_i(\cdot)\}$ and $\{\gamma_i\}$. We list some of their properties here:

$\mathrm{(a)}$ For any $i\in\mathbb{N}$, the elements in the family $\mathcal{F}_i=\{\overline{B}(x,\gamma_i):x\in K_i\}$ are disjoint due to the choice of $\gamma_i.$ Each element in  $\mathcal{F}_{i+1}$ is a subset of $\overline{B}(x,\frac{\gamma_i}{2})$ for some $x\in K_i$. To see this, first recall that $K_{i+1}=\bigcup\limits_{x\in K_{i}}E_{i+1}(x)$, $E_{i+1}\subset K\cap B(x,\frac{\gamma_i}{4})$ and $0<\gamma_{i+1}<\frac{\gamma_i}{4}$. Then for every $x_0\in K_{i+1}$, there exists a point $x\in K_i$ such that $x_0\in E_{i+1}(x)\subset K\cap \overline{B}(x,\frac{\gamma_i}{4})$. And hence $\overline{B}(x_0,\gamma_{i+1})\subset\overline{B}(x_0,\frac{\gamma_i}{4})\subset \overline{B}(x,\frac{\gamma_i}{2}).$


$\mathrm{(b)}$ For each $x\in K_i$ and $z\in\overline{B}(x,\gamma_i)$, 
\begin{align}\label{eq-B}
  \overline{B}(z,m_i(x),\varepsilon,\phi)\cap\bigcup\limits_{y\in K_i\backslash\{x\}}\overline{B}(y,\gamma_i)=\emptyset
\end{align}
and
\begin{align}\label{eq-muB}
  \mu_i(\overline{B}(x,\gamma_i))=e^{-m_i(x)s}\leq\sum\limits_{y\in E_{i+1}(x)}e^{-m_{i+1}(y)s}\leq (1+2^{-i-1})\mu_i(\overline{B}(x,\gamma_{i})),
\end{align}
 where $E_{i+1}(x)=B(x,\gamma_i)\cap K_{i+1}.$

By (a), for every $j\geq i,\;K_j\subset\bigcup\limits_{x\in K_i}\overline{B}(x,\frac{\gamma_i}{2})$ which implies $\mu_i(F_i)\leq\mu_{j}(F_i),\;\forall F_i\in\mathcal{F}_i.$ 
By \eqref{eq-B}, 
$$\mu_i(F_i)\leq\mu_{i+1}(F_i)=\sum\limits_{F\in\mathcal{F}_{i+1}:F\subset F_i}\mu_{i+1}(F)\leq(1+2^{-i-1})\mu_i(F_i),\;F_i\in\mathcal{F}_i.$$

Applying the above inequalities repeatedly, it holds that for any $j>i,$
\begin{align}\label{eq-mu}
\mu_i(F_i)\leq\mu_j(F_i)\leq\prod\limits_{n=i+1}^j(1+2^{-n})\mu_i(F_i)\leq C\mu_i(F_i),\;\forall F_i\in\mathcal{F}_i.
\end{align}
where $C=\prod\limits_{n=1}^\infty(1+2^{-n})<\infty.$

Let $\tilde{\mu}$ be a limit point of $\{\mu_i\}$ in the weak-star topology and $\tilde{K}=\bigcap\limits_{n=1}^\infty\overline{\bigcup\limits_{i\geq n}K_i}$. It's easy to check that $\tilde{\mu}$ is supported on $\tilde{K}$.  Moreover, $\tilde{K}$ is a compact subset of $K$ because 
each $\overline{\bigcup\limits_{i\geq n}K_i}$ is a closed subset of the compact set $K$. 

By \eqref{eq-mu}, for each $x\in K_i$, we have 
$$e^{-m_i(x)s}=\mu_i(\overline{B}(x,\gamma_i))\leq\widetilde{\mu}(\overline{B}(x,\gamma_i))\leq C\mu_i\left(\overline{B}(x,\gamma_i)\right)=Ce^{-m_i(x)s}.$$
Particularly, we have $$1<\sum_{x\in K_1}e^{-m_1(x)s}=\sum_{x\in K_1}\mu_1(B(x,\gamma_1))\leq\mu_1(\tilde{K})\leq\tilde{\mu}(\tilde{K})\leq\sum\limits_{x\in K_1}C\mu_1(B(x,\gamma_1))\leq2C.$$
Note that $\tilde{K}\subset \bigcup_{x\in K_i}\overline B(x, \frac{\gamma_i}{2})$. 
By \eqref{eq-B}, for each $x\in K_i$ and $z\in\overline{B}(x,\gamma_i)$, we have $$\tilde{\mu}(\overline{B}(z,m_i(x),\varepsilon,\phi))\leq\tilde{\mu}(\overline{B}(x,\frac{\gamma_i}{2}))\leq Ce^{-m_i(x)s}.$$
For each $z\in \tilde{K}$ and $i\in\N$, $z\in B(x,\frac{\gamma_i}{2})$ for some $x\in K_i$. Hence 
$$\tilde{\mu}(\overline{B}(z,m_i(x),\varepsilon,\phi))\leq Ce^{-m_i(x)s}.$$

Define $\mu=\tilde{\mu}/{\tilde{\mu}(\tilde{K})}$. Then $\mu\in M(\tilde{K})\subset M(K)\subset M(X)$. For each $z\in \tilde{K},$ there exists a sequence $t_i\nearrow\infty$ such that $\mu(B(z,t_i,\varepsilon,\phi))\leq Ce^{-t_is}/\tilde{\mu}(\tilde{K})$. It follows that $$\overline{h}_{\mu}(\phi)=\int\lim_{\varepsilon\rightarrow0}\limsup_{t\rightarrow +\infty}-\frac{1}{t}\log\mu(B(x,t,\varepsilon,\phi))d\mu\geq s.$$
\end{proof}
Combining Lemma \ref{lemma4-1} and Lemma \ref{lemma4-2}, we complete the proof of Theorem \ref{main}.

{\bf Acknowledgements}
This work is supported by the National University Innovation Contest Grant of Nanjing University under the guidance of Professor Dou Dou. The authors thank Professor Fei Yang for his encouragement.

\end{document}